\documentclass[11pt]{article}
\usepackage{latexsym}
\usepackage{amssymb}
\usepackage{amsmath}
\usepackage{mathrsfs}
\usepackage{amsfonts}
\usepackage{amsthm}

\newcommand{\equald}{\mbox{$ \;\stackrel{\cal D}{=}\; $}}

\newcommand{\RL}{{\mathbb R}}
\newcommand{\NN}{{\mathbb N}}
\newcommand{\IN}{{\mathbb Z}}

\newcommand{\IND}{{\mathbb I}}

\newcommand{\PR}{\mbox{\rm Pr}} 
\newcommand{\Ent}{\mbox{Ent}}

\newcommand{\la}{\lambda}

\newcommand{\CPolQ}{\mbox{\rm CP}(\la,Q)} 
\newcommand{\CPm}{{\rm CP}_{\la,Q}} 
 
\newcommand{\Pmj}{{\rm P}_{\la Q_j}}

\newcommand{\smallsum}{{\textstyle \sum}}

\def\be{\begin{eqnarray}}
\def\ee{\end{eqnarray}}
\def\ben{\begin{eqnarray*}}
\def\een{\end{eqnarray*}}

\def\elabel#1{\label{e:#1}}

\def\sq{$\Box$}

\def\qed{\ifmmode\sq\else{\unskip\nobreak\hfil
\penalty50\hskip1em\null\nobreak\hfil\sq
\parfillskip=0pt\finalhyphendemerits=0\endgraf}\fi\par\medbreak}

\newsavebox{\junk}
\savebox{\junk}[1.6mm]{\hbox{$|\!|\!|$}}

\def\til={{\widetilde =}}

\def\half{{\mathchoice{\textstyle \frac{1}{2}}%
{\frac{1}{2}}%
{\hbox{\tiny $\frac{1}{2}$}}%
{\hbox{\tiny $\frac{1}{2}$}} }}

 \def\eq#1/{(\ref{#1})}

\def\eq#1/{(\ref{e:#1})}

\newcommand{\beqn}[1]{\notes{#1}%
\begin{eqnarray} \elabel{#1}}

\newcommand{\eeqn}{\end{eqnarray} }

\newcommand{\beq}[1]{\notes{#1}%
\begin{equation}\elabel{#1}}

\newcommand{\eeq}{\end{equation}} 

\def\bdes{\begin{description}}
\def\edes{\end{description}}

\newcommand{\levy}{\mbox{L\'evy}} 
 
\def\notes#1{}

\title{MEASURE CONCENTRATION\\ FOR COMPOUND POISSON DISTRIBUTIONS}

\author{I.\ KONTOYIANNIS
\footnote{Supported in part
by a Sloan Foundation Research Fellowship
and by NSF grant \#0073378-CCR.}\\
{\em Division of Applied Mathematics, Brown University,}\\
{\em 182 George Street, Providence, RI 02912, USA}\\
Email: {\tt yiannis@dam.brown.edu}
\and 
M.\ MADIMAN\\
{\em Division of Applied Mathematics, Brown University,}\\
{\em 182 George Street, Providence, RI 02912, USA.}\\
Email: \texttt{mokshay@dam.brown.edu}
}

\date{}

\begin{document}
\bibliographystyle{plain}
\maketitle



\noindent
{\bf Abstract: } 
We give a simple development of the concentration 
properties of compound Poisson measures 
on the nonnegative integers. A new modification of the 
Herbst argument is applied to an appropriate modified 
logarithmic-Sobolev inequality to derive new
concentration bounds. When the measure of interest 
does not have finite exponential moments, these bounds
exhibit optimal {\em polynomial} decay.
Simple new proofs are also given for earlier results 
of Houdr{\'e} (2002) and Wu (2000).

\medskip


\medskip



\section{Introduction}

Concentration of measure is a well-studied phenomenon,
and in the past 30 years or so it has been explored
through a wide array of tools and techniques;
\cite{talagrand:95}\cite{ledoux:book}\cite{mcdiarmid:98}
offer broad introductions.
Results in this area are equally well
motivated by theoretical questions 
(in areas such as geometry, functional analysis
and probability),
as by numerous applications in 
different fields including the analysis of algorithms,
mathematical physics and empirical processes in statistics.

From the probabilistic point of view, measure concentration
describes situations where a random variable is 
strongly concentrated around a particular value. This 
is typically quantified by the rate of decay of the probability 
that the random variable deviates from that value
(usually its mean or median) by a certain amount. 
As a simple concrete example consider a function
$f(W)$ of a Poisson$(\la)$ random variable $W$;
if $f:\IN_+\to\RL$ is 1-Lipschitz, i.e.,
$|f(k)-f(k+1)|\leq 1$ for all $k\in\IN_+=\{0,1,2,\ldots\}$,
then \cite{bobkov-ledoux:98},
\be
\PR\big\{f(W)-E[f(W)]>t\big\}
\leq\exp\Big\{-\frac{t}{4}\log\Big(1+\frac{t}{2\la}
\Big)\Big\}.
\label{eq:BLcom}
\ee
Although the distribution of 
$f(W)$ may be quite complex,
(\ref{eq:BLcom}) provides a simple, explicit 
bound on the probability that it deviates from 
its mean by an amount $t$. This is a general theme:
Under appropriate conditions, it is possible to 
derive useful, accurate bounds of this type
for a large class of random variables with complex 
and often only partially known distributions. 
We also note that the consideration of Lipschitz 
functions is motivated by applications,
but it is also related to fundamental concentration 
properties captured by isoperimetric inequalities
\cite{ledoux:book}.  

The bound (\ref{eq:BLcom}) was established 
in \cite{bobkov-ledoux:98}
using the so-called ``entropy method,'' pioneered by
Ledoux \cite{ledoux:96}\cite{ledoux:97}\cite{ledoux:book}.
The entropy method consists of two steps.
First, a (possibly modified) logarithmic-Sobolev inequality 
is established 
for the distribution of interest. Recall that, for an 
arbitrary probability measure $\mu$ and any nonnegative
function $f$ on the same space, 
the entropy functional $\Ent_\mu(f)$ 
is defined by 
$$\Ent_\mu(f)=
\smallint f\log f d\mu-(\smallint f d\mu) \log(\smallint f d\mu),$$
whenever all the above integrals exist.
In the case of the Poisson, Bobkov and Ledoux 
\cite{bobkov-ledoux:98}
proved the following modified log-Sobolev inequality:
Writing ${\rm P}_\la$ for the Poisson$(\la)$
measure, for any function $f:\IN_+\to\RL$
with positive values,
$$
\Ent_{{\rm P}_\la}(f)\leq\la E_{{\rm P}_\la}\Big[\frac{1}{f}|Df|^2\Big],
$$
where $Df(k)=f(k+1)-f(k),$ $k\geq 0$, is the discrete gradient,
and $E_\mu$ denotes the expectation operator with 
respect to a measure $\mu$.
In fact, they also established the following
sharper bound which we will use below; for any
function $f$ on $\IN_+$,
\be 
\Ent_{{\rm P}_\la} \big(e^f\big) \leq \la E_{{\rm P}_\la} \Big[
e^f 
\big\{ 
|Df|e^{|Df|}-e^{|Df|}+1 
\big\}
\Big].
\label{eq:poi-lsi-2}
\ee

The second step in the entropy method is the so-called
Herbst argument: Starting from some Lipschitz function $f$, 
the idea is to use the modified log-Sobolev inequality
to obtain an upper bound on the entropy of
$e^{\tau f}$, and from that to 
deduce a differential inequality for the 
moment-generating function $G(\tau)=E[e^{\tau f}]$
of $f$. Then, solving the differential inequality 
yields an upper bound on $G(\tau)$,
and this leads to a concentration bound via 
Markov's inequality.

Our main goal in this work is to carry out a similar
program for an arbitrary {\em compound Poisson}
measure on $\IN_+$. Recall that for any $\la>0$ 
and any probability measure $Q$ on the natural 
numbers $\NN=\{1,2,\ldots\}$, the compound Poisson 
distribution $\CPolQ$ is the distribution of the random 
sum 
\ben
Z\equald \sum_{i=1}^W X_i,
\een
where $W\sim$ Poisson$(\la)$ and the $X_i$ 
are independent random variables with distribution
$Q$ on $\NN$, also independent of $W$;
we denote the $\CPolQ$ measure by $\CPm$.
The class of compound Poisson distributions is 
much richer than the one-dimensional Poisson
family. In particular, 
the $\CPolQ$ law inherits its tail behavior 
from $Q$: $\CPolQ$ has
finite variance iff $Q$ does, it has exponentially
decaying tails iff $Q$ does, and so on \cite{sato:book}.
It is in part from this versatility of tail behavior 
that the compound Poisson distribution draws its 
importance in many applications. 
Alternatively,
$\CPolQ$ is characterized
as the infinite divisible law
without a Gaussian component
and with $\levy$ measure $\la Q$.

From the above discussion we observe that
the Herbst argument is heavily dependent on the
use of moment-generating functions, a fact which
implicitly assumes the existence of exponential 
moments. Our main contribution 
is a modification of the Herbst argument for the
case when the random variables of interest do
{\em not} satisfy such exponential integrability 
conditions. We derive what appear to be perhaps
the first concentration inequalities for a class
of infinitely divisible random variables that
have finite variance but do not have finite
exponential moments.
Apart from the derivation of the
present results, the modified Herbst argument 
is applicable in a variety of other cases
and may be of independent interest.
In particular, this approach can be applied 
to prove dimension-free inequalities 
for compound Poisson vectors, as well as 
power-law concentration bounds for more general 
infinitely divisible laws.

Our starting point is the following modified log-Sobolev inequality
for the compound Poisson measure $\CPm$.

\medskip

\noindent
{\bf Theorem 1. } [{\sc Modified Log-Sobolev Inequality for
Compound Poisson Measures}]$\;$
For any $\la>0$, any probability measure $Q$ on $\NN$
and any bounded $f:\IN_+\to\RL$,
\be
\Ent_{\CPm} \big(e^f\big) \leq 
\la \sum_{j\geq 1} Q_j E_{\CPm}\Big[
e^f 
\big\{ 
|D^jf|e^{|D^jf|}-e^{|D^jf|}+1 
\big\}
\Big],
\label{eq:cpoi-lsi}
\ee
where $D^jf(k)=f(k+j)-f(k)$, for $j,k\in\IN_+$.

\medskip

This can be derived easily from \cite[Cor~4.2]{wu:00b} of Wu,
which was established using elaborate stochastic calculus
techniques. In Section~3 we also give an alternative,
elementary proof, by tensorizing the Bobkov-Ledoux result
(\ref{eq:poi-lsi-2}). Note the elegant similarity between the 
bounds in (\ref{eq:poi-lsi-2}) and (\ref{eq:cpoi-lsi}).

We then apply our modified Herbst argument
to establish concentration bounds for 
$\CPolQ$ measures under various assumptions on the
tail behavior of $Q$. These
are stated in Section~2 and proved in Section~4.
For example, we establish the following polynomial
concentration result. Recall that a function $f:\IN_+\to\RL$
is $K$-Lipschitz if $|f(j+1)-f(j)|\leq K$ for all $j\in\IN_+$.

\medskip

\noindent
{\bf Corollary 2.} [{\sc Polynomial Concentration}]
Suppose that $Z$ has $\CPolQ$ distribution where $Q$ 
has finite moments up to order $L$,
$$L=\sup\big\{\tau\geq 1\,:\,\smallsum_{j\geq 1}j^\tau\,Q_j<\infty\big\}>1,$$
and write $q_r$ for its integer moments,
$$q_r=\smallsum_{j\geq 1}j^r\,Q_j,\;\;\;\;r\in\NN.$$
If $f:\IN_+\to\RL$ is $K$-Lipschitz,
then for any positive integer $n<L$ and any $t>0$ we have,
$$\PR\big\{|f(Z)-E[f(Z)]|>t\big\}\leq A \cdot B^n\cdot t^{-n},
$$
where for the constants $A,B$ we can take,
\ben
A
&=&
\exp\Big\{\la\sum_{r=1}^n\binom{n}{r}q_rK^r-\la n\log K\Big\}\\
B
&=&
2|f(0)|+2K \la q_1 + 1.
\een

Various stronger and more general results are given in Section~2. 
There, at the price of more complex constants, 
we get bounds which, for large $t$, are of (the optimal)
order $t^{-L+\delta}$ 
for any $\delta>0$. Moreover, since the only property 
of the compound Poisson distribution used in the proof 
is that it satisfies the functional inequality of Theorem~1, 
similar bounds are immediately seen to hold for any
measure that satisfies such an inequality.
Note that although the bound of
Corollary~2 is not useful for small $t$, it is 
in general impossible to obtain meaningful results 
for arbitrary $t>0$. For example, if $f$ 
is the identity function and $Z\sim$ Poisson$(\la)$ 
where $\la$ is
of the form $m+1/2$ for an integer $m$, 
then $|Z-E(Z)|\geq 1/2$ with probability 1;
a more detailed discussion is given in Section~2.

As noted above, these appear to be some of 
the first {\em non-exponential} concentration 
bounds that have been derived, with the few recent 
exceptions discussed next. 
Of the extensive current literature on concentration,
our results are most closely related to the work of Houdr{\'e} 
and his co-authors. Using sophisticated technical 
tools derived from the ``covariance representations''
developed in \cite{houdre-P-S:98}\cite{houdre-privault:02}, 
Houdr{\'e} \cite{houdre:02} obtained concentration bounds 
for Lipschitz functions of infinitely divisible 
random vectors with finite {\em exponential} moments.
In \cite{houdre-marchal:04}, truncation and explicit 
computations were used to extend these results to
the class of stable laws on $\RL^d$, and the 
preprint \cite{breton-houdre-P:pre} extends them further 
to a large class of functionals on Poisson space.
To our knowledge, the results in 
\cite{houdre-marchal:04}\cite{breton-houdre-P:pre} 
are the only concentration bounds with power-law 
decay to date. But when specialized to scalar random 
variables they only apply to distributions with 
infinite variance, whereas our results hold for 
compound Poisson random variables with a finite 
$L$th moment for any $L>1$.
Although the methods of 
\cite{houdre-marchal:04}\cite{breton-houdre-P:pre} 
as well as the form of the results themselves 
are very different from those derived here, 
some more detailed comparisons are possible 
as outlined in Section~2.
Finally, the recent paper \cite{BBLM:05} 
contains a different extension of the Herbst argument 
to certain situations where exponential moments 
do not exist. The focus there is on
{\em moment} inequalities for functions 
of independent random variables, primarily 
motivated by statistical applications.

\section{Concentration Bounds}

The following result is the main motivation for this paper.
It illustrates the potential for using the Herbst argument 
even in cases where the existence of exponential moments fails
or cannot be assumed.

\medskip

\noindent
{\bf Theorem 3.} [{\sc Power-law Concentration}] 
Suppose that $Z$ has $\CPolQ$ distribution where $Q$ has 
finite moments up to order $L$,
$$L=\sup\big\{\tau\geq 1\,:\,\smallsum_{j\geq 1}j^\tau\,Q_j<\infty\big\}>1,$$
and write 
$q_1=\smallsum_{j\geq 1}j\,Q_j$
for its first moment.
\begin{itemize}
\item[{\bf (i)}]
If $f:\IN_+\to\RL$ is $K$-Lipschitz,
then for any $t>0$ and $\epsilon>0$ we have,
\be
\Pr\big\{|f(Z)-Ef(Z)|> t\big\} \leq 
\exp
\bigg\{
\inf_{0<\alpha<L}
\Big[
I_{\epsilon}(\alpha)+\alpha
\log\Big(\frac{2|f(0)|+2K \la q_1 +\epsilon}{t}\Big)
\Big] 
\bigg\},
\label{eq:pol-res2}
\ee
where 
\ben
I_{\epsilon}(\alpha)
&=&
\la\sum_{j\geq 1} Q_j \big\{ C_{j,\epsilon}^{\alpha}  
-1-\alpha \log C_{j,\epsilon} \big\}\\
C_{j,\epsilon}
&=&
1+\frac{jK}{\epsilon}.
\een
\item[{\bf (ii)}]
The upper bound (\ref{eq:pol-res2}) is meaningful 
(less than 1) iff
$t>T:=2|f(0)|+2K \la q_1+\epsilon$, 
and then,
\ben
\Pr\big\{|f(Z)-Ef(Z)|> t\big\} 
\leq
\exp\bigg\{ - 
\int_0^{\log(t/T)} 
i_{\epsilon}^{-1}(s) ds  \bigg\}
\een
where
$i_{\epsilon}(\alpha) := I_{\epsilon}'(\alpha)
=
\la \sum_{j\geq 1} Q_j 
[C_{j,\epsilon}^\alpha-1]\log C_{j,\epsilon}.$
\end{itemize}

\noindent
{\bf Remarks. }

1. Taking $\alpha=L-\delta$ for any $\delta>0$
in the exponent of (\ref{eq:pol-res2}), we get a
bound on the tails
of $f(Z)$ of order $t^{-(L-\delta)}$ for 
large $t$. 
By considering 
the case where $f$ is the identity function
$f(k)=k$, $k\in\IN_+$,
we see that this power-law behavior is in fact
optimal. In particular, this shows that the
tail of the $\CPolQ$ law decays like the tail of $Q$, 
giving a quantitative version of
a classical result from \cite{sato:book}.

2. As will become evident from the proof, Theorem~3 holds 
for any random variable $Z$ with 
law $\mu$ instead of $\CPm$, as long as $\mu$ satisfies 
the log-Sobolev inequality of Theorem~1
with respect to some probability measure $Q$ on $\NN$ 
and some $\la>0$, and assuming that $\mu$ has finite 
moments up to order $L$.
The bound (\ref{eq:pol-res2}) remains 
exactly the same, except that the first moment 
$M_1=E[Z]$ of $\mu$ replaces $\la q_1$.

3. Integrability properties follow immediately from the theorem: 
For any $K$-Lipschitz function $f$,
$E_{\CPm}[|f|^{\tau}]<\infty$ for all $\tau<L$, and the same
holds for any law $\mu$ as in the previous remark.

\medskip

Since the support of $\CPolQ$ is $\IN_+$, we would
naturally expect the range of $f$ to be highly 
disconnected. Therefore, to somewhat simplify the expression 
in the exponent of (\ref{eq:pol-res2}) next we concentrate
on the (typical) class of functions $f:\IN_+\to\RL$ whose 
mean under $\CPm$ is not in the range of $f$:

\medskip

\noindent
{\bf Corollary 4.} [{\sc Power-law Concentration for Nice $f$}] 
Suppose that $Z$ has $\CPolQ$ distribution where $Q$ has 
finite moments up to order $L>1,$
and write $q_1$ for its first moment.
If $f:\IN_+\to\RL$ is $K$-Lipschitz
and there exists
$\epsilon>0$ such that
$$|f(j)-E[f(Z)]|\geq\epsilon,\;\;\;\;\mbox{for all }j\in\IN_+,$$
then for any $t>0$ we have,
\be
\Pr\big\{|f(Z)-Ef(Z)|> t\big\} \leq 
\exp
\bigg\{
\inf_{0<\alpha<L}
\Big[
I_{\epsilon}(\alpha)+\alpha
\log(D/t)
\Big] 
\bigg\},
\label{eq:nice-f}
\ee
where $I_{\epsilon}(\alpha)$ is defined as in Theorem~3, and
$D:=E|f(Z)-E[f(Z)]|$.

\medskip

\noindent
{\bf Remarks.}

4. Similarly to Theorem~3, this corollary gives quantitative 
bounds on the tail of $f(Z)$ of the order of $t^{-(L-\delta)}$
for any $\delta>0$. Also, the same result holds for any
law $\mu$ as in Remark~2.

5. 
The exponent in (\ref{eq:nice-f}) becomes negative
exactly when $t>D$, for the 
same reasons as in Theorem~2. On the other hand, it is
obvious that any bound can only be 
useful for $t>D_0:=\min_{k\in\IN_+} |f(k)-E[f(Z)]|$, 
since the probability that $|f(Z)-E[f(Z)]|\geq D_0$ 
is equal to one. Moreover, $D$ and $D_0$ 
coincide in many special cases, as, e.g., when 
the range of $f$ is a lattice in $\RL$ and its mean 
$E[f(Z)]$ is on the midpoint between two lattice points.
In this sense, the restriction $t>D$ is quite natural.

6. The expression $2|f(0)|+2K\la q_1$ in Theorem~3 
is simply an upper bound to the constant $D=E|f(Z)-E[f(Z)]|$ 
appearing in Corollary~4. In both cases, when $L>2$ it is
possible to obtain potentially sharper results by 
bounding $D$ above using Jensen's inequality by,
$$ \big[ \{K^2\la q_2+\{|f(0)|+K\la q_1\}^2 \big]^{\half},$$
where $q_2$ is the second moment of $Q$. Similar
expressions can be derived in the case of higher moments.

7. The most closely related results to our power-law 
concentration bounds appear to be in the recent 
preprint \cite{breton-houdre-P:pre}.\footnote{The
results in \cite{breton-houdre-P:pre} are stated 
in the much more general setting of functionals 
on an abstract Poisson space.
Using the Wiener-Ito decomposition, any 
infinitely divisible random variable can be 
represented as a Poisson stochastic integral, 
which in turn can be realized as a ``nice'' 
functional on Poisson space.}
The relevant bounds in \cite{breton-houdre-P:pre} 
specialized to Lipschitz functions of 
$\CPolQ$ random variables require that the
probability measure $Q$ be {\em non-atomic}, 
which excludes all the cases we consider. 
But shortly after the first writing of this 
paper C.~Houdr{\'e} in a personal communication
informed us that this assumption can be removed 
by an appropriate construction. 
The details have not been checked by us, 
but in the following comparison we assume 
that it does not change the statements 
in \cite{breton-houdre-P:pre}.
The main assumptions in \cite{breton-houdre-P:pre}
are that the random variable of interest has infinite
variance, and also certain growth conditions.
Because of the infinite-variance assumption,
the majority of the results in this paper
(corresponding to $L>2$) apply to cases that
are not covered in \cite{breton-houdre-P:pre}.
As for the growth conditions, they are convenient 
to check in several important special classes, e.g., 
for $\alpha$-stable laws on $\RL$, but they 
can be unwieldy in the compound Poisson case, 
especially as they depend on $Q$ in an intricate way. 
On the other hand, if $Q$ has infinite variance,
\cite[Cor.~5.3]{breton-houdre-P:pre} gives 
optimal-order bounds, including the case
when $Q$ has infinite mean, for which our results 
do not apply.

\medskip

Next we show how the Herbst argument can be used to recover
precisely a result of \cite{houdre:02} in the case 
when we have exponential moments.

\medskip

\noindent
{\bf Theorem 5.} [{\sc Exponential Concentration}] \cite{houdre:02}
Suppose that $Z$ has $\CPolQ$ distribution where $Q$ has 
finite exponential moments up to order $M$,
$$M=\sup\big\{\tau\geq 0\,:
\,\smallsum_{j\geq 1}e^{\tau j}\,Q_j<\infty\big\}>0.$$
If $f:\IN_+\to\RL$ is $K$-Lipschitz, 
then for any $t>0$ we have,
\be
\Pr\big\{f(Z)-Ef(Z)> t\big\} 
\;\leq \;
\exp
\Big\{
\inf_{0<\alpha<M/K}
[
H(\alpha)-\alpha t
] 
\Big\}
\;=\;
\exp
\Big\{
-\int_0^t h^{-1}(s)\,ds
\Big\},
\label{eq:expthm}
\ee
where $H(\alpha):=\la \sum_{j\geq 1}Q_j[e^{\alpha Kj}-1-\alpha K j]$,
and $h^{-1}$ is the inverse of 
$h(\alpha):= H'(\alpha).$

\medskip

\noindent
{\bf Remarks.}

8. Theorem~1 of \cite{houdre:02} gives concentration
bounds for a class of infinitely divisible laws with finite
exponential moments, and in the compound Poisson case it reduces
precisely to (\ref{eq:expthm}), which 
also applies to any random variable $Z$ 
whose law satisfies the result of Theorem~1. 
It is also interesting to note that Theorem~5 
can be derived by applying
\cite[Prop~3.2]{wu:00b} to a compound Poisson random 
variable (constructed via the Wiener-Ito decomposition),
and then using Markov's inequality.


9. Theorems~3 and~5 
easily generalize to 
H{\"o}lder continuous
functions.
In the discrete setting of $\IN_+$, 
$f$ is $K$-Lipschitz iff
it is H{\"o}lder continuous for every 
exponent $\beta\geq 1$ with the same 
constant $K$. But if $f$ is H{\"o}lder continuous 
with exponent $\beta< 1$, this more stringent 
requirement makes it possible to strengthen
Theorem~3 and Theorem~5,
by respectively redefining,
$C_{j,\epsilon}=
1+\frac{j^\beta K}{\epsilon}$, and
\ben
H(\alpha)
=
\la \sum_{j\geq 1}Q_j
\big[e^{\alpha K j^\beta}-1-\alpha K j^\beta\big].
\een

10. While all our power-law results dealt with two-sided
deviations, the bound in Theorem~5 
is one-sided. The reason for this
discrepancy is that 
the last step in all the relevant proofs is an 
application of Markov's inequality, which leads
us to restrict attention to nonnegative random 
variables. When exponential moments exist, the 
natural consideration of the exponential
of the random variable takes care of this issue,
but in the case of regular moments we are forced
to take absolute values.

\section{Proof of Theorem~1}

An alternative representation
for the law of a $\CPolQ$ random variable $Z$ is in terms of the series
\be\label{eq:represent}
Z\equald \sum_{j=1}^\infty j\,Y_j,\;\;\;\;\;\;
Y_j\sim\mbox{Poisson}(\la Q_j),
\ee
where the $Y_j$ are independent.

For each $n$, let $\mu_n$ denote the joint 
(product) distribution
of $(Y_1,\ldots,Y_n)$.
In this instance,
the tensorization property 
of the entropy 
\cite{bobkov:96}\cite{ledoux:97}\cite{ledoux:book}
can be expressed as
\be
\mbox{Ent}_{\mu_n}(G)\leq\sum_{j=1}^n
E\Big[
\mbox{Ent}_{\Pmj}
\big(G_j(Y_1,\ldots,Y_{j-1},\cdot,Y_{j+1},\ldots,Y_n)\big)
\Big],
\label{eq:tensorize2}
\ee
where $G:\IN_+^n\to\RL_+$
is an arbitrary function,
and the entropy on the right-hand side 
is applied to the restriction 
$G_j$ of $G$ to its $j$th co-ordinate.
Now given an $f$ as in the statement of the
theorem, define the functions
$G:\IN_+^n\to\RL_+$ and $H:\IN_+^n\to\RL_+$
by
$$H(y_1,\ldots,y_n)=f\Big(\sum_{k=1}^n ky_k\Big),
\;\;\;y_1^n\in\IN_+^n,$$
and $G=e^H$. Let 
$\bar{\mu}_n$ denote the distribution
of the sum $S_n:=\sum_{k=1}^nkY_k$ and
write $H_j:\IN_+\to\RL$ for
the restriction of $H$ to the
variable $y_j$ with the remaining
$y_i$'s fixed. 
Applying (\ref{eq:tensorize2}) to $G$ we obtain,
\ben
\mbox{Ent}_{\bar{\mu}_n}(e^f) 
\;=\;
\mbox{Ent}_{\mu_n}(G)
&\leq& 
\sum_{j=1}^n
E\Big[
\mbox{Ent}_{\Pmj}
\Big( G_j  (Y_1,\ldots,Y_{j-1},\cdot,Y_{j+1},\ldots,Y_n) \Big)
\Big]\\
&=& 
\sum_{j=1}^n
E\Big[
\mbox{Ent}_{\Pmj}
\big(
e^{H_j(Y_1,\ldots,Y_{j-1},\cdot,Y_{j+1},\ldots,Y_n)}
\big)
\Big].
\een
Using 
the Bobkov-Ledoux inequality (\ref{eq:poi-lsi-2})
to bound each term in the above sum,
and noting that, trivially,
$DH_j(y_1,\ldots,y_n)=D^j f(\sum_{k=1}^nky_k),$
\be
\mbox{Ent}_{\bar{\mu}_n}(e^f)
&\leq&
\sum_{j=1}^n \la Q_j
E_{\mu_n}\Big[
e^{H}
\big\{ 
|DH_j|e^{|DH_j|}-e^{|DH_j|}+1 
\big\}
\Big]
\nonumber\\
&=&
\la\sum_{j=1}^n Q_j
E_{\bar{\mu}_n}\Big[
e^{f}
\big\{ 
|D^jf|e^{|D^jf|}-e^{|D^jf|}+1 
\big\}
\Big],
\nonumber\\
&\leq&
\la\sum_{j=1}^\infty Q_j
E_{\bar{\mu}_n}\Big[
e^{f}
\big\{ 
|D^jf|e^{|D^jf|}-e^{|D^jf|}+1 
\big\}
\Big],
\label{eq:finite-n}
\ee
where the last inequality follows from the fact 
that $xe^x-e^x+1\geq 0$ for $x\geq 0$.

Finally, we want to take the limit as $n\to\infty$ in
(\ref{eq:finite-n}). 
Since $\bar{\mu}_n\Rightarrow\CPm$ 
as $n\to\infty$ by (\ref{eq:represent}), and since
$f$ is bounded, by bounded convergence
\be
\mbox{Ent}_{\bar{\mu}_n}(e^f)
\to
\mbox{Ent}_{\CPm}(e^f),
\;\;\;\;n\to\infty.
\label{eq:limit1}
\ee
Similarly, changing the order of summation and expectation
in the right-hand side of (\ref{eq:finite-n}) by Fubini,
taking $n\to\infty$ by bounded convergence,
and interchanging the order again, 
it converges to 
\ben
\la\sum_{j\geq 1} Q_j
E_{\CPm}\Big[
e^{f}
\big\{ 
|D^jf|e^{|D^jf|}-e^{|D^jf|}+1 
\big\}
\Big].
\een
This together with (\ref{eq:limit1}) 
implies that (\ref{eq:finite-n}) yields the required result upon
taking $n\to\infty$.
\qed

\section{Concentration Proofs}

For notational convenience we define the function
$\eta(x):=xe^x -e^x +1,$
$x\in\RL,$
and note that it is non-negative; 
it achieves its minimum at 0; 
it is strictly convex on $(-1,\infty)$ and 
strictly concave on $(-\infty,-1)$; 
it decreases from 1 to 0 as $x$ increases
to zero, and it is
increasing to infinity for $x>0$.

The main technical ingredient of the paper is 
the following proposition, which is based
on a modification of the Herbst argument.

\medskip

\noindent
{\bf Proposition 7.}
Suppose that $Z$ has $\CPolQ$ distribution where $Q$ has 
finite moments up to order $L>1$.
If $f:\IN_+\to\RL$ is bounded and $K$-Lipschitz,
then for $t>0,$ $\epsilon>0$ and $\alpha\in(0,L)$, we have,
\ben
\Pr\{|f(Z)-Ef(Z)|> t\} \leq 
\exp
\Big\{
I_{\epsilon}(\alpha)+\alpha E[\log g_\epsilon(Z)]
-\alpha\log t
\Big\},
\een
where $I_{\epsilon}(\alpha)$ is defined as in Theorem~3 and
\ben
g_{\epsilon}(x):=
 |f(x)-E[f(Z)]| \,\IND_{\{ \,|f(x)-E[f(Z)]|\geq\epsilon \,\} } 
+ \epsilon \, \IND_{ \{ \,|f(x)-E[f(Z)]|<\epsilon \,\} }.
\een


\noindent
{\sc Proof of Proposition~7. }
Since $f$ is bounded, by its definition $g_\epsilon$ is also
bounded above by $2\|f\|_\infty+\epsilon$ and below by $\epsilon$.
Therefore, the moment generating function 
$G(\tau):=E[g_\epsilon(Z)^\tau]$
is well-defined for all $\tau> 0$. Moreover, since both
$g_\epsilon$ and $\log g_\epsilon$ are bounded, dominated
convergence justifies the following differentiation under
the integral,
\ben
G'(\tau)=E\bigg[\frac{\partial}{\partial\tau} 
e^{\tau \log g_{\epsilon}(Z)}\bigg] 
= E \big[ g_{\epsilon}(Z)^\tau \log g_{\epsilon}(Z)\big]
\een
so we can relate $G(\tau)$ to the entropy of $g_\epsilon^\tau$,
\be
\Ent_{\CPm} (g_{\epsilon}^{\tau})= \tau G'(\tau)-G(\tau)\log G(\tau)
= \tau^2 G(\tau)\frac{d}{d\tau} \bigg[ \frac{\log G(\tau)}{\tau} \bigg] .
\label{eq:expand}
\ee

In order to bound this entropy we will apply Theorem~1 to the 
function $\phi(x):=\tau\log g_\epsilon(x)$. First we observe
that $g_\epsilon$ can be written as the composition 
$g_\epsilon=h\circ(f-E[f(Z)])$, where it is easy to
verify that the function 
$h(x):=|x|\IND_{\{|x|\geq\epsilon\}}
+\epsilon\IND_{\{|x|<\epsilon\}}$ is
1-Lipschitz. And since $f$ is $K$-Lipschitz
by assumption, $g_\epsilon$
is itself $K$-Lipschitz. Hence we can 
bound $D^j\phi$ as
\ben
D^j \phi(x)
=
\tau \log\bigg|\frac{g_{\epsilon}(x+j)}{g_{\epsilon}(x)}\bigg| 
\leq
\tau \log\bigg( 1+ \bigg|\frac{D^j g_{\epsilon}(x)}
{g_{\epsilon}(x)}\bigg|\bigg) 
\leq
\tau \log\bigg( 1+ \frac{jK}{\epsilon}\bigg) 
=
\tau \log C_{j,\epsilon}.
\een
The same argument also yields a corresponding lower bound,
so that $|D^j \phi(x)| \leq \tau \log C_{j,\epsilon}.$
Applying Theorem~1 to $\phi$ gives,
\ben
\Ent_{\CPm} (g_{\epsilon}^{\tau})
= 
\Ent_{\CPm} (e^\phi)
\leq
\la\sum_{j\geq 1}Q_jE_{\CPm}[e^\phi\eta(|D^j\phi|)]
\leq
\la G(\tau)\sum_{j\geq 1}Q_j
\eta(\tau\log C_{j,\epsilon}),
\een
since $\eta(x)$ is increasing for $x\geq 0$.
Combining this with (\ref{eq:expand})
we obtain the following 
differential inequality valid for all $\tau>0$:
\ben
\frac{d}{d\tau} \bigg[ \frac{\log G(\tau)}{\tau} \bigg] 
\leq \la  \sum_{j\geq 1} Q_j \frac{\eta(\tau \log C_{j,\epsilon})}{\tau^2}.
\een
To solve,
we integrate with respect to $\tau$ 
on $(0,\alpha]$ to obtain, for any $\alpha<L$,
\ben
\frac{\log G(\alpha)}{\alpha} - E[\log  g_{\epsilon}(Z)]
&\leq&
\la \sum_{j\geq 1} Q_j \int_{0}^{\alpha} 
\frac{\eta(\tau \log C_{j,\epsilon})}{\tau^2} d\tau  \\
&=&
\la \sum_{j\geq 1} Q_j \log C_{j,\epsilon} 
\int_{0}^{\alpha\log C_{j,\epsilon}} \frac{\eta(s)}{s^2}\, ds  \\
&=&  
\la \sum_{j\geq 1} Q_j \log C_{j,\epsilon} 
\bigg[ \frac{e^s -1-s}{s} \bigg]_{0}^{\alpha\log C_{j,\epsilon}} \\
&=&  
I_\epsilon(\alpha)/\alpha,
\een
or, equivalently,
\be
G(\alpha) \leq \exp \big\{
\alpha E[\log  g_{\epsilon}(Z)] + I_{\epsilon}(\alpha) \big\},
\label{eq:MGF}
\ee
where the exchange of sum and integral is
justified by Fubini's theorem since all the
quantities involved are nonnegative.
To complete the proof we observe that 
$g_{\epsilon}\geq |f-E[f(Z)]|$, so that
by (\ref{eq:MGF}) and an application of Markov's inequality, 
\ben 
\hspace{0.95in}
\PR\big\{|f(Z)-E[f(Z)]|>t\big\}
&\leq&
\PR\big\{g_\epsilon(Z)>t\big\}\\
&=&
\PR\big\{g_\epsilon(Z)^\alpha>t^\alpha\big\}\\
&\leq&
t^{-\alpha}\cdot G(\alpha)\\
&\leq&
\exp \bigg\{ 
I_{\epsilon}(\alpha) 
+\alpha E[\log  g_{\epsilon}(Z)] 
-\alpha \log t 
\bigg\}.
\hspace{0.95in}
\Box
\een
 
\medskip

Using Proposition~7 we can prove our main
results, Theorem~3 and Corollaries~2 and~4.

\medskip

\noindent
{\sc Proof of Theorem~3. }
The first step is to bring 
the upper bound in Proposition~7 into a more tractable form.
Observe that by its definition, 
$g_{\epsilon}(x)\leq |f(x)-E[f(Z)]| +\epsilon,$
so that, by Jensen's inequality,
for a function $f$ satisfying the hypotheses 
of Proposition~7,
\be
E[\log  g_{\epsilon}(Z)] \leq \log E[g_{\epsilon}(Z)]
\leq
\log\Big[
E\big\{|f(Z)-E[f(Z)]|\big\} +\epsilon 
\Big].
\label{eq:g-bound}
\ee
Thus the upper bound in Proposition~7 can be weakened to
\be
\PR\big\{
|f(Z)-E[f(Z)]|>t
\big\}
\leq
\exp \bigg\{  I_{\epsilon}(\alpha) 
+ \alpha\log \big( D +\epsilon \big) 
-\alpha \log t \bigg\} ,
\label{eq:Delta}
\ee
where $D:=E\big\{|f(Z)-E[f(Z)]|\big\}.$
Next we use the Lipschitz property of $f$ to obtain an upper 
bound for the above exponent 
which is uniform over all $f$ with $f(0)$ 
fixed. Since $f(j)\in[f(0)-Kj,f(0)+Kj]$, we have
$|f(j)|\leq|f(0)|+Kj$, and hence
\ben
D
\leq
2E|f(Z)|
\leq 
2|f(0)|+2K\la q_1,
\een
where we used the fact that the mean of the $\CPolQ$ law 
is $\la q_1$. Substituting in (\ref{eq:Delta})
and taking the infimum over 
$\alpha$ yields the required result (\ref{eq:pol-res2}),
and it only remains to remove the boundedness assumption on $f$.
But since the bound itself only depends on $f$ via
$f(0)$ and $K$, truncating $f$
at level $\pm n$ and passing to the limit $n\to\infty$
proves part (a).

With $T=2|f(0)|+2K \la q_1 +\epsilon,$ in order 
to evaluate the exponent
\be
\inf_{0<\alpha<L}
\big[
I_{\epsilon}(\alpha)+\alpha
\log\big(T/t\big)
\big]
\label{eq:exp-tmp}
\ee
in (\ref{eq:pol-res2}), we calculate the 
first two derivatives of 
$I_\epsilon(\alpha)$ with respect to $\alpha$ as,
\ben
I'_\epsilon(\alpha)
=
\la\sum_{j\geq 1} Q_j[C_{j,\epsilon}^\alpha-1]\log C_{j,\epsilon}
\;\;\;\;\;
\mbox{and}
\;\;\;\;\;
I''_\epsilon(\alpha)
=
\la\sum_{j\geq 1} Q_jC_{j,\epsilon}^\alpha(\log C_{j,\epsilon})^2,
\een
where the exchange of differentiation and expectation 
is justified by dominated convergence;
observe that, since $C_{j,\epsilon}>1$, both are
positive for all $\alpha>0$. In particular, since
$I_\epsilon(\alpha)>0$, 
the exponent (\ref{eq:exp-tmp}) can only be negative
(equivalently, the bound in (\ref{eq:pol-res2}) can
only be less than 1) if the second term in (\ref{eq:exp-tmp})
is negative, i.e., if $t>T$. On the other hand, since 
$I'_\epsilon(0)=0$
and 
$I''_\epsilon(\alpha)>0$ for all $\alpha$, we see that
$I_\epsilon(\alpha)$ is locally quadratic around $\alpha=0$.
This means that, as long as $t>T$, choosing $\alpha$
sufficiently small we can make (\ref{eq:exp-tmp}) negative,
therefore the bound of the theorem is meaningful precisely
when $t>T$.

To obtain the alternative representation,
fix any $\epsilon>0$ and set
$i_{\epsilon}(\alpha)=I_{\epsilon}'(\alpha)$.
Since $I_\epsilon''(\alpha)$ is strictly positive,
for $t>T$ the expression
$I_{\epsilon}(\alpha)+\alpha\log(T/t)$ is 
uniquely minimized at $\alpha^*>0$ which solves
$i_{\epsilon}(\alpha)=\log(t/T)>0$. 
Hence, for all $t>T$, integrating by parts,
\ben
\min_{0<\alpha<L} \Big[ I_{\epsilon}(\alpha)+\alpha\log(T/t) \Big] 
&=& I_{\epsilon} (\alpha^*) + \alpha^* \log(T/t)\\
&=& \int_0^{\alpha^*} i_{\epsilon}(s)ds  + \alpha^* \log(T/t)\\
&=& \int_0^{i_\epsilon(\alpha^*)}  
x \,di_{\epsilon}^{-1}(x)  + \alpha^* \log(T/t)\\
&=& 
i_\epsilon(\alpha^*)i_\epsilon^{-1}\big(i_\epsilon(\alpha^*)\big)
- \int_0^{i_\epsilon(\alpha^*)} 
i_{\epsilon}^{-1}(x) dx   + \alpha^* \log(T/t) \\
&=&  - \int_0^{\log(t/T)} i_{\epsilon}^{-1}(x) dx ,
\een
which proves part (b).
\qed

\medskip

\noindent
{\sc Proof of Corollary~4. }
The proof is identical to that
of Theorem~3, with the only difference
that, since here we simply have
$g_{\epsilon}(x)=  |f(x)-E[f(Z)]|$ for all $x$,
we can replace the bound 
(\ref{eq:g-bound}) by
$ E[\log g_{\epsilon}(Z)] \leq \log D$,
where $D=E|f(Z)-E[f(Z)]|$.
Proceeding
as before
gives the result.
\qed

\medskip

\noindent
{\sc Proof of Corollary 2. }
This is an application of
Theorem~3 for specific values of 
$\alpha$ and $\epsilon$: 
Bounding the infimum by the 
value at $\alpha=n$ and taking
$\epsilon=1$,
\be
\Pr\big\{|f(Z)-Ef(Z)|> t\big\} \leq 
\exp
\bigg\{
I_{1}(n)+n
\log\Big(\frac{2|f(0)|+2K \la q_1 +1}{t}\Big)
\Big] 
\bigg\}.
\label{eq:consequence}
\ee
Using the binomial theorem to expand
$I_1(n)$,
\ben
I_1(n)\;=\;
\la\sum_{j\geq 1} Q_j  \big\{ C_{j,1}^n - 1 -n\log C_{j,1}\big\}
&=&
\la\sum_{j\geq 1} Q_j \big\{ (1+jK)^n - 1\big\} 
-\la n\sum_{j\geq 1} Q_j\log (1+jK)\\
&\leq&
\la\sum_{j\geq 1} Q_j 
\sum_{r=1}^n\binom{n}{r}
(jK)^r 
-\la n\sum_{j\geq 1} Q_j[\log j+\log K]\\
&\leq &
\la\sum_{r=1}^n
\binom{n}{r}
K^r q_r
-\la n\log K.
\een
Substituting this bound into (\ref{eq:consequence})
and rearranging yields the result.
\qed

\medskip

Next we go on to prove the exponential concentration
result Theorem~5 using the classical Herbst argument
in conjunction with the modified log-Sobolev inequality 
of Theorem~1.

\medskip

\noindent
{\sc Proof of Theorem 5. } We proceed similarly to 
the proof of Proposition~7.
Assume $f$ is a bounded and $K$-Lipschitz,
and let $F(\tau)=E[\exp\{\tau f(Z)\}]$,
$\tau>0$ be the moment-generating function of $f(Z)$.
Dominated convergence justifies 
the differentiation
\ben
\label{exp-diff}
F'(\tau)=E[f(Z)e^{\tau f(Z)}],
\een
so we can relate $F'(\tau)$ to the entropy of $e^{\tau f}$ by
\be
\Ent_{\CPm}(e^{\tau f}) 
=
\tau F'(\tau)-F(\tau)\log F(\tau) 
= 
\tau^2 F(\tau) \frac{d}{d\tau} \bigg[ \frac{\log F(\tau)}
{\tau} \bigg].
\label{eq:relate}
\ee
Since $f$ is $K$-Lipschitz, the function 
$g:=\tau f$ is $\tau K$-Lipschitz, so that 
$|D^jg|\leq\tau Kj$. 
Applying Theorem~1 to $g$,
\ben
\Ent_{\CPm}(e^{\tau f}) 
=
\Ent_{\CPm}(e^{g}) 
\leq
\la\sum_{j\geq 1} Q_j E\Big[e^{g(Z)}\eta(|D^jg(Z)|)\Big]
\leq
\la F(\tau)\sum_{j\geq 1} Q_j \,\eta(j\tau K).
\een
Combining this with (\ref{eq:relate}) yields
\ben
\frac{d}{d\tau}\bigg[\frac{\log F(\tau)}{\tau} \bigg]
\leq \la   \sum_{j\geq 1} Q_j \bigg\{ 
\frac{\eta(j\tau K)}{\tau^2} \bigg\} ,
\een
and integrating with respect to $\tau$ from 0 to $\alpha>0$
we obtain
\ben
\frac{\log F(\alpha)}{\alpha} - E[f(Z)]
&\leq&
\la\int_0^{\alpha} \sum_{j\geq 1} Q_j \bigg\{ 
\frac{\eta(j\tau K)}{\tau^2} \bigg\} d\tau \\
&=&
\la\sum_{j\geq 1} jK Q_j \int_0^{j\alpha K} 
\frac{\eta(s)}{s^2}\;ds  \\
&=&
\la\sum_{j\geq 1} Q_j 
\bigg[ \frac{e^{jK\alpha}-1-jK\alpha}{\alpha} \bigg] ,
\een
where the exchange of the sum and integral is justified
by Fubini's theorem since the integrand is nonnegative.
Therefore, we have the following a bound on the moment-generating 
function $F$,
\be
F(\alpha) \leq \exp\big\{\alpha E[f(Z)] + H(\alpha)\big\},
\;\;\;\;\alpha>0,
\ee
where $H(\alpha)=\la  \sum_j Q_j \big[ e^{jK\alpha}-1-jK\alpha \big] $.
An application of Markov's inequality now gives
\ben
\PR\big\{ 
f(Z)-E[f(Z)] >t \big\} 
&\leq&
e^{-\alpha t}
E\Big[\exp\big\{\alpha [f(Z)-E[f(Z)])]\big\}\Big]\\
&=& 
e^{-\alpha t}
F(\alpha)
e^{-\alpha E[f(Z)]}\\
&\leq& 
\exp\big\{ H(\alpha) -\alpha t \big\} .
\een

The removal of the boundedness assumption
is a routine truncation argument as 
in the proof of Theorem~3 or in
\cite{bobkov-ledoux:98}\cite{houdre:02}.
In order to obtain the best bound for the deviation
probability, we minimize the exponent
over $\alpha\in(0,M/K)$. This yields the first expression
in Theorem~5; the second representation
follows from a standard 
argument as in the last part of the proof 
of Theorem~3 or \cite{houdre:02}. 
\qed


\section*{Acknowledgments}
We thank Christian Houdr{\'e} for sending us a copy 
of the preprint \cite{breton-houdre-P:pre},
and we also thank Christian Houdr{\'e} and Andrew Barron
for numerous comments on these and other related results.



\end{document}